\documentclass{article}
\usepackage[T1]{fontenc}
\usepackage{amsmath}
\usepackage{enumitem}
\usepackage{graphicx}
\usepackage{array}
\usepackage{booktabs}
\usepackage{graphicx}
\usepackage{changepage}
\usepackage{longtable}
\usepackage[skip=0.1\baselineskip]{caption}

\begin{document}

\title{Collatz Conjecture: Patterns Within}
\author{H. Nelson Crooks, Jr. \\
		Chigozie Nwoke}
\date{May 2022}
\maketitle
\begin{abstract}
Collatz Conjecture sequences increase and decrease in seemingly random fashion. By identifying and analyzing the \emph{forms} of numbers, we discover that Collatz sequences are governed by very specific, well-defined rules, which we call \emph{cascades}.
\end{abstract}

\newpage

\tableofcontents
\newpage
\listoftables
\listoffigures
\newpage

\section{Introduction}
\subsection{Background and Approach} 
The Collatz Conjecture was posed by Lothar Collatz \text{c. 1937}. The Conjecture asks whether the iterative algorithm

\[
C_{i+1}=
\begin{cases}
3*C_{i}+1,	& \text{for } C_{i} \text{ odd} \\
C_{i}/2, & \text {for } C_{i} \text{ even}
\end{cases}
\]
converges to 1 for all natural starting numbers $ C_{0} $. 

\bigskip
The Collatz Conjecture has been checked by computer and found to reach 1 for all numbers $C_{0} \leq 2^{68}\approx 2.95 x 10^{20}$ [1], but has not been proven to be true for all natural numbers.

In this paper, we develop a unique method of describing natural numbers and analyzing their behavior under the Collatz algorithm. This methodology allows us to draw conclusions about large groups of numbers. We seek to understand the structure of the random increase and decrease of a Collatz sequence.

We will identify a number of interesting patterns in Collatz sequences we have found during our analysis. 

\subsection{Hypothesis}
The sequence for every number iterates through a value that is smaller than the initial value. The iteration after which a value smaller than the starting number is reached for the first time is called the \emph{stopping time}. If we can prove that all natural numbers have a finite stopping time, we can prove by induction that the Collatz Conjecture is true.

We analyze Collatz sequences for groups of starting numbers using forms, cascades, and columns.

\section{Number Forms}
\subsection{Forms}
This analysis of the Collatz Conjecture is based on the concept of number \emph{forms}, which consist of a \emph{base} (a natural number) times an \emph{index} (a whole number) plus an \emph{offset} (a whole number). We define

\begin{table}[!htbp]
\begin{tabular}{l l l}
\emph{standard form} as 		& $C=2^{p}n+2^{p-1}-1$; 		& $p \in N,n\in W$\\
\emph{non-standard form} as & $C=2^{p}n+2^{p}-1$; 		& $p \in N,n \in W$\\
\emph{composite form} as 	& $C=2^{p}n+f$; 				& $p \in N, n \in W, f \in W,$\\
  								&									& $f <2^{p},f \neq 2^{p-1}-1,f \neq 2^{p}-1$  \\
\emph{mixed form} as & $ C=kn+f$; & $k \in N, k\neq \text{a power of 2}, n \in W,$\\
& &$ f \in W, f<k.$\\

\end{tabular}
\end{table}

\newpage
\textbf{Lemma 1:} All natural numbers can be expressed in standard form by the equation:

\begin{table}[!htbp]
\centering 		
\begin{tabular}{l l r}
$C=2^{p} n+2^{p-1}-1;$ 		& $p \in N, n \in W (N \cup \{0\})$  		&  (1)\\ 
\end{tabular}
\end{table}
	
This equation was developed independently by the authors, and was previously identified by Mehendale [5] and Cadogan [2], and also describes even numbers when p=1.
\bigskip
See Appendix I for proof of Lemma 1.

The standard form of a natural number is determined by the value of p in equation (1). Small values of p and corresponding bases, offsets, and standard forms are shown in Table \ref{table:StdForms}.

\begin{table}[!htbp]
\caption{Standard Number Forms} 			
\centering 										
\begin{tabular}{c c c c} 						
\hline\hline									
\rule{0pt}{4mm}
p 	& Base	&Offset	& Standard Form	\\		
\hline 											
1 	& 2 	& 0 	& 2n 				\\ 				
2 	& 4 	& 1 	& 4n+1			\\ 					
3 	& 8 	& 3 	& 8n+3			\\ 					
4 	& 16 	& 7 	& 16n+7			\\ 					
5 	& 32 	& 15 	& 32n+15			\\						
\hline 											
\end{tabular}
\label{table:StdForms}
\end{table}

Cadogan identified standard form numbers in Table 1 of [2] with\newline  $A_{1}=4n+1, A_{2}=8n+3, etc.$  

We define a \emph{\#-form} as a number for which the $base = \#$, e.g., a 4-form number is of the form $4n+1$.

We call the sequence of form numbers (the bases) for natural numbers (4, 2, 8, 2, 4, 2, 16, 2, 4, 2, 8, 2, 4, 2, 32, …) a \emph{form pattern}. We’ll see this form pattern again later in our analysis. 

See Appendix II for a method to determine the standard form of a number.

\bigskip
\textbf{Lemma 2:} Every natural number can be described in standard form by a unique combination of p and n. Conversely, the combination of p and n for each natural number is unique. See Appendix III for proof of Lemma 2.

\subsection{Collatz Steps}
The Collatz algorithm requires that if a number is odd we multiply the number by 3 and add 1. We call this an \emph{odd step}. The result of an odd step is always an even natural number of the mixed form $3t+1, t \in \text{\emph{odd }} N$.  We note that the result of an odd step cannot be a multiple of 3.

If a number is even, we divide the number by two. We call this an \emph{even step}. The result of an even step can be any natural number, odd or even.

Because multiples of 3 cannot result from an odd step, multiples of 3 in a Collatz sequence can only be the first number(s) of the sequence. 
\subsection{Collatz Cycles}
Since an odd step must always be followed by an even step, we now define an \emph{odd cycle} as an odd step followed by an even step: multiply by 3 and add 1, then divide by 2. We define an \emph{even cycle} as an even step. A \emph{general cycle} is an odd or even cycle as appropriate.

We further define a \emph{\#-cycle} as a cycle that starts with a number of the form \emph{\#}. For example, a 4-cycle starts with a number that is a 4-form number, such as $13 (= 4*3 +1)$. 

All odd cycles increase the value of a number by slightly more than half. A 2-cycle (the only even cycle) reduces the value of a number by half (we call it a 2-cycle to distinguish it from the even step at the end of all odd cycles).

\section{Cascades}
\subsection{Analysis of a general cycle}
We now analyze what happens to an odd number of the standard form \newline $C=2^{p} n+2^{p-1}-1$ when a general cycle is applied. $C_{i}, C_{i+1}, C_{i+2},\dotsc$ represent successive steps in the Collatz algorithm.

\begin{table} [!htbp]
\caption{Analysis of a General Cycle} 			
\centering 										

\begin{tabular}{l l l l} 						
\hline\hline

\rule{0pt}{4mm}
		& $C=2^{p} n+2^{p-1}-1$ 								& $\text{for odd } C, p \in N,$		&  (1) \\
		&															& $p \geq 2, n \in W$ 					&\\ 			
Step 1	& $C_{i} \text{ is odd by definition here } \dotsc$ 		& so apply an odd step 				& \\
 		& $C_{i+1}=3C_{i}+1$ 									& 										& \\
 		& $C_{i+1}=3(2^{p}n+2^{p-1}-1)+1$ 					& Substitute (1) for $C_{i}$ 			& \\
 		& $C_{i+1}=3(2^{p}n)+3(2^{p-1})-3(1)+1$ 				& Distribute the 3 						& \\
 		& $C_{i+1}=3(2^{p}n)+3(2^{p-1})-2$ 					& Combine like terms 				& \\
 		& $C_{i+1}=2(3(2^{p-1}n)+3(2^{p-2})-1)$ 				& Factor out 2: 						& (2) \\
 		& 															& $\text{after 1 iteration }=C_{i+1}$ 	& \\
Step 2	& $C_{i+1} \text{ is even, }\dotsc$ 						& so apply an even step 				& \\
 		& $C_{i+2}=C_{i+1}/2$ 									& 										& \\
 		& $C_{i+2}=2(3(2^{p-1}n)+3(2^{p-2})-1)/2$ 				& Substitute (2) for $C_{i+1}$ 	& \\
 		& $C_{i+2}=3(2^{p-1}n)+3(2^{p-2})-1$ 				& Cancel the 2's 						& \\
 		& Manipulating the result $\dotsc$ 						& 										& \\
 		& $C_{i+2}=3(2^{p-1}n)+2(2^{p-2})+1(2^{p-2})-1$ 	& Split $3(2^{p-2})$  into 			& \\
 		&															& $2(2^{p-2}) + 1(2^{p-2})$			& \\ 
 		& $C_{i+2}=3(2^{p-1}n)+1(2^{p-1})+1(2^{p-2})-1$ 	& Move the 2 into the 					& \\
 		&															& exponent 							& \\
 		& $C_{i+2}=2^{p-1}(3n+1)+2^{p-2}-1$ 				& Combine like terms; 					& (3) \\
 		& 															& after 2 iterations $=C_{i+2}$ 		& \\
\hline
\end{tabular}
\label{table:GenCycle}
\end{table}

\newpage
\textbf{Applying a general cycle to an odd number results in a number of the next lower form (4-form is lower than 8-form) with the index of the initial number multiplied by 3 and increased by 1!} 
(Cadogan identified this as Theorem 1 [2]). See Appendix IV for an example. 

\subsection{Cascades}
Because we evaluated a general cycle, it is clear that a subsequent cycle will reduce the form and increase the value of the result of the previous cycle. This will continue until a 2-form number is reached (an even number) at which time a 2-cycle will divide the number by 2. We call this pattern of successive cycles from the starting form through and including the 2-cycle a \emph{cascade}. We further identify a \emph{\#-cascade} as a cascade that starts with a \#-form number. 

Once a cascade starts, the sequence is rigidly defined by alternating odd and even steps until a 2-form number is reached, at which time another even step is required. Each cascade is comprised of $2p-1$ steps (p-1 odd steps and p even steps), where p is the power of 2 in the base of the starting number. We see that the base of the number at each step of the cascade is even through the 4-cycle of the cascade and odd after the final 2-cycle. Thus, the odd/even parity of the result of the cascade is indeterminate and depends on the parity of n. Based on values of n, half of all cascades result in an even number and half in an odd number. The result of a cascade with an odd value of p has the same parity as n. The result of a cascade with an even value of p has the opposite parity of n. 

\vspace{0.5cm}
\textbf{The seemingly random behavior of a Collatz sequence can now be seen as successive odd cascades interspersed with stand-alone 2 cascades.} (We call them \emph{2-cascades}, even though they involve only an even step, to distinguish them from the \emph{2-cycle} at the end of each higher-form cascade.) 

The distinctive sawtooth increase/decrease portions of a plot of a Collatz sequence are the visual representations of cascades. \textbf{The specific form (values of p and n) of the number that results from a cascade determines whether the sequence will then go up or down, and by how much. This, we believe, is the key to understanding the Collatz Conjecture.}

\subsection{Cascade Transforms}
Since the steps in a cascade are rigidly constrained, we can shortcut our sequence calculations by substituting the result of a \#-cascade for the calculations of the individual cycles within the cascade. Table \ref{table:CascadeTransforms} shows the cascade transforms for some small standard forms.

\begin{table}[!htbp]
\caption{Cascade Transforms} 			
\centering 										
\begin{tabular}{c c c c} 						
\hline\hline 									
\rule{0pt}{4mm}
p 	& Base & Standard Form	& Transform	\\			
\hline 											
1  	& 2 	& 2n 				& n  			\\ 
2 	& 4 	& 4n+1			& 3n+1		\\ 					
3 	& 8 	& 8n+3			& 9n+4		\\ 					
4	& 16 	&16n+7			& 27n+13		\\ 				
5 	& 32 	& 32n+15			& 81n+40		\\				
6 	& 64 	& 64n+31			& 243n+121	\\ 					
7 	& 128	&128n+63			&729n+364	\\ 						
\hline 											
\end{tabular}
\label{table:CascadeTransforms}
\end{table}

\vspace{12pt}
Note the pattern in these cascade transforms; 

$2^{p}n + 2^{p-1} - 1$ transforms to $3^{p-1}n + (3^{p-1} - 1)/2.$ 

\vspace{12pt}
Cascade transforms are always mixed-form numbers, where the base is not a power of 2.

\newpage
\subsection{Form Patterns in Cascade Transforms}
When we analyze the form pattern of cascade transforms for 8-cascades we see a familiar pattern.
\vspace{-5mm}
\begin{table}[!htbp]
\caption{Forms of 8-Cascade Transforms} 			
\centering 										
\begin{tabular}{c c c c c} 						
\hline\hline 									
\rule{0pt}{4mm}
		& Input      	& Transform	&                 	& 				\\
         	& Standard	& Mixed       	& Standard   	& Standard 	\\			
Index 	& Form       	& Form        	& Form of     	& Base of 		\\			
 n  		& 8n+3      	& 9n+4        	& Transform 	& Transform	\\			
\hline 											
0    	&   3  			&    4 			& 2(2)         	& 2  			\\ 
1    	& 11   			&   13 			& 4(3)+1     	& 4  			\\ 
2    	& 19   			&   22 			& 2(11)       	& 2				\\ 
3    	& 27   			&   31 			& 64(0)+31 	& 64			\\ 
4    	& 35   			&   40 			& 2(20)       	& 2				\\ 
5    	& 43   			&   49 			& 4(12)+1  	& 4				\\
6    	& 51   			&   58 			& 2(29)      	& 2				\\ 
7    	& 59   			&   67 			& 8(8)+3    	& 8				\\
8    	& 67   			&   76 			& 2(38)      	& 2   			\\ 
9    	& 75   			&   85 			& 4(21)+1  	& 4   			\\ 
10  	& 83   			&   94 			& 2(47)      	& 2   			\\
11  	& 91   			& 103 			& 16(6)+7  	& 16   			\\ 
12 		& 99   			& 112 			& 2(56)      	& 2   			\\ 
13 		& 107 			& 121 			& 4(30)+1  	& 4				\\ 
14 		& 115 			& 130 			& 2(65)      	& 2  			\\ 
15 		& 123 			& 139 			& 8(17)+3  	& 8   			\\ 
16 		& 131 			& 148 			& 2(74)      	& 2				\\ 
17 		& 139 			& 157 			& 4(39)+1  	& 4				\\ 
\hline 											
\end{tabular}
\label{table:8CascadeTransforms}
\end{table}

The form pattern of the 8-transforms (2, 4, 2, 64, 2, 4, 2, 8, 2, 4, 2, 16, 2, 4, 2, 8, 2, 4,…) looks like the form pattern for standard forms of sequential numbers, shifted by some amount. The form patterns of cascade transforms match the form patterns of natural numbers for small standard forms when the index is shifted by values shown in Table \ref{table:PatternShifts}. This behavior appears to apply to all 8- and higher transforms.

\begin{table}[!htbp]
\caption{Matching Form Patterns in Cascade Transforms} 			
\centering 										
\begin{tabular}{ c c c c c c} 						
\hline\hline 									
\rule{0pt}{4mm}
In     		& From		& To      		& Cascade   		& From   	&To  			\\
Form 		& Index \#	& Index \# 	& Transform 	 	& Index \# & Index \# 	\\
\hline										

4n+1		& 0 		& 524,286 		& 3n+1 			& 174,763 	& 699,049  	\\ 
8n+3 		& 0 		& 524,286 		& 9n+4 			& 407,780 	& 932,066 		\\
16n+7 		& 0 		& 524,286 		& 27n+13 			& 135,927 	& 660,213 		\\
32n+15 	& 0 		& 779,928 		& 81n+40 			& 220,072 	& 1,000,000 	\\
64n+31 	& 0 		& 577,117 		& 243n+121 		& 422,883 	& 1,000,000 	\\
128n+63 	& 0 		& 524,286 		& 729n+364 		& 315,724 	& 840,010 		\\
\hline
\end{tabular}
\label{table:PatternShifts}
\end{table}
\newpage
\subsection{Maximum Cascade Starts (MCS)}
Every natural number can be the result of a cascade, and every natural number has a \emph{maximum cascade start} (maximum form, minimum value) that transforms to that number. Trümper calls them \emph{Collatz Backward Series} in [7].

We remember that each odd cycle in a cascade results in a number of the next-lower-form with the index of the previous number multiplied by 3 and increased by 1.

We find the maximum cascade start (MCS) for a number by working backwards from the number. We call this reverse cascade a \emph{ladder}.

\begin{table}[!htbp]
\caption{Finding Maximum Cascade Start} 			
\centering 										
\begin{tabular}{l p{10cm}} 														
\hline\hline
Step 1		& Every cascade ends with a 2-cycle, so we start by multiplying the number by 2. \\ 
 			& \\
 Step 2		& If the index of this number is of the mixed form $3t+1, t \in W$, then there is a higher 
form number that will cycle to this number. We use t from the index in a form 1 higher than the form of the number to generate the next-higher-form number in the ladder.\\
 			&\\
 Step 3		& Repeat Step 2 until the index is not of the mixed form  $3t+1, t \in W$, then go to Step 4. \\
			&  \\
 Step 4		& If the index of this number is NOT of the mixed form  $3t+1, t \in W$, then there is no higher-form number that will cycle to this number, and we have found the maximum cascade start.\\
\hline
\end{tabular}
\label{table:FindMCS}
\end{table}

Since the cascade structure is rigidly defined, we can use this method to find all odd values in the ladder. 

We don’t need to find the even numbers in the odd cycles of the reverse cascade.

Just as there is only one path for each cascade, there is only one path for each ladder.

\newpage
The example below helps understand the process to find the maximum cascade start for a number.

\begin{table}[!htbp]
\caption{Finding Maximum Cascade Start Example}
\centering 										
\begin{tabular}{l c c c c} 														
\hline\hline 									
\rule{0pt}{4mm}
		& Form & Number 	& Index 	& Next \\
\hline
\rule{0pt}{4mm}
		& $64(0)+31=$	&$31$		&					&\\
		&&&&\\
Step 1 	& $2(31)+ 0=$		& $62$ 	& $31=3(10)+1$ 	&up to next-\\
		&&&&higher form\\
		&&&&\\
Step 2	& \multicolumn{4}{l}{The index of 62 is of the form $3t+1$, so there is a higher-form}\\ 
		&\multicolumn{4}{l}{number that will cycle to 31.}\\
		&&&&\\
Step 2 	& \multicolumn{4}{l}{Using 10 as the index for the next number in the ladder results in}\\
		&&&&\\
 		& $4(10)+ 1=$ 	& $41$ 	& $10=3(3)+1$ 	&up to next-\\
		&&&&higher form\\
		&&&&\\
Step 2 	& $8(3)+ 3=$ 		& $27$ 	& $3=3(1)+0$ 		& end of\\
		&&&&reverse cascade\\
		&&&&\\
Step 3 & \multicolumn{4}{l}{The index of 27 is not of the form $3t+1$, so there is no higher-form}\\ 					  		&\multicolumn{4}{l}{number that will cycle to 27.}\\
		&&&&\\
Step 4 & \multicolumn{4}{l}{Thus, 27 is the maximum cascade start of 31.}\\ 
		&&&&\\
		&\multicolumn{4}{l}{Checking that the 27 cascade does lead to 31:}\\
		&\multicolumn{4}{l}{$27=8(3)+3$}\\
		&\multicolumn{4}{l}{$82=2(41)$}\\
		&\multicolumn{4}{l}{$41=4(10)+1$}\\
		&\multicolumn{4}{l}{$124=2(62)$}\\
		&\multicolumn{4}{l}{$62=2(31)$}\\
		&\multicolumn{4}{l}{$31=64(0)+31$ (end of cascade)}\\
\hline
\end{tabular}
\label{table:FindMCSExample}
\end{table}

\noindent Maximum cascade start values for some small numbers are shown in Table \ref{table:MCSValues}.

\begin{table}[!htbp]
\caption{Maximum Cascade Starts for Some Small Numbers} 			
\centering 										
\begin{tabular}{ccc} 						
\hline\hline 									
\rule{0pt}{4mm}
						&Maximum Cascade		&Standard Base of\\
Cascade Ending Value	&Starting Value			&Maximum Cascade\\
(Mixed Form)			&(Max Form, Min Value)	&Starting Value\\
\hline										
\rule{0pt}{4mm}
$28 = 3(9) + 1 $ 		& $37 = 4(9) + 1$ 		& $	4$	\\
$29 = 1(29)$ 			& $	58 = 2(29)	$ 			& $2$	\\
$30 = 1(30)$ 			& $	60 = 2(30)$ 			& $	2$	\\
$31 = 9(3) + 4	$ 		& $ 27 = 8(3) + 3$ 		& $	8$	\\
$32 = 1(32)$ 			& $	64 = 2(32)	$ 			& $2$	\\
$33 = 1(33)$ 			& $	66 = 2(33)$ 			& $	2$	\\
$34 = 3(11) + 1$ 		& $	45 = 4(11) + 1$ 		& $	4$	\\
$35 = 1(35)$ 			& $	70 = 2(35)	$ 			& $2$	\\
$36 = 1(36)$ 			& $	72 = 2(36)$ 			& $	2$	\\
$37 = 3(12) + 1$ 		& $	49 = 4(12) + 1$ 		& $	4$	\\
$38 = 1(38)$ 			& $	76 = 2(38)$ 			& $	2$	\\
$39 = 1(39)$ 			& $	78 = 2(39)$ 			& $	2$	\\
$40 = 81(0) + 40$		& $	15 = 32(0) + 15$ 		& $	32$\\
$41 = 1(41)$ 			& $	82 = 2(41)$ 			& $	2$	\\
$42 = 1(42)$ 			& $	84 = 2(42)$ 			& $	2$	\\
$43 = 3(14) + 1$ 		& $	57 = 4(14) + 1$ 		& $	4$	\\
\hline 											
\end{tabular}
\label{table:MCSValues}
\end{table}

\newpage
We can see a modified form pattern in these maximum values (4, 2, 2, 8, 2, 2, 4, 2, 2, 4, 2, 2, 32, 2, 2, 4,…). An interesting characteristic of the maximum cascade start values is that the index of each is NOT of the form 3n+1. If it were in that form there would be a higher cascade start. 

\subsection{Primary Maximum Cascade Starts (PMCS)}
If successive maximum cascade starts (MCS) are applied to any number that is not a multiple of 3 (find the MCS of the MCS, etc.) it appears that an odd multiple of 3 is eventually reached. Proof of this is offered by Trümper in [7]. For example, successive MCS of 28 are: 28, 37, 49, 43, 57. This appears to be true for any natural number that is not a multiple of 3 and has been verified by the authors for all numbers up to 1,700,000.

Because a multiple of 3 can only be the initial value(s) of a cascade, we call the cascade that starts with an odd multiple of 3 the “primary maximum cascade start” for a given number.

\begin{table}[!htbp]
\caption{Primary Maximum Cascade Starts for Some Small Numbers} 			
\centering 										
\begin{tabular}{ccc} 						
\hline\hline 									
\rule{0pt}{4mm}
						&							&Standard Base of\\
Cascade Ending Value	&Primary Maximum		&Maximum Cascade\\
(Mixed Form)			&Cascade Starting Value	&Starting Value\\
\hline										
\rule{0pt}{4mm}
$28 = 3(9) + 1	$ 		& $57 = 4(14) + 1$ 		& $	4$	\\
$29 = 1(29)$			& $	51 = 8(6)+3	$ 		& $8$	\\
$30 = 1(30)$ 			& see note 				& see note	\\
$31 = 9(3) + 4	$ 		& $27 = 8(3) + 3$ 		& $	8$	\\
$32 = 1(32)$ 			& $	75 = 8(9)+3	$ 		& $8$	\\
$33 = 1(33)$ 			& see note 					& see note	\\
$34 = 3(11) + 1$ 		& $	45 = 4(11) + 1$ 		& $	4$	\\
$35 = 1(35)$ 			& $	93 = 4(23)+1	$ 		& $4$	\\
$36 = 1(36)$ 			& see note 					& 	see note	\\
$37 = 3(12) + 1$ 		& $	57 = 4(14) + 1$ 		& $	4$	\\
$38 = 1(38)$ 			& $	39 = 16(2)+7$ 		& $	16$	\\
$39 = 1(39)$ 			& see note					& 	see note	\\
$40 = 81(0) + 40$ 	& $	15 = 32(0) + 15$ 		& $	32$	\\
$41 = 1(41)$ 			& $	171 = 8(21)+3$ 		& $	8$	\\
$42 = 1(42)$ 			& see note					& 	see note	\\
$43 = 3(14) + 1$ 		& $	57 = 4(14) + 1$ 		& $	4$	\\
\multicolumn{3}{c}{Note: Multiples of 3 cannot result from an odd cascade.}\\
\hline 											
\end{tabular}
\label{table:PMCSValues}
\end{table}
\newpage
\section{Seeds}
The only way for the Collatz algorithm to reach 1 is through a power of 2, which goes to 1 via repeated 2 cycles. We call the number that leads to a power of 2 a \emph{seed}. The first few seeds and the power of 2 they lead to are shown below.

\begin{table}[!htbp]
\caption{Seeds (Small Values)}
\centering 										
\begin{tabular}{p{15mm}ll}					
\hline\hline
\rule{0pt}{4mm}
&$1 \Rightarrow 4 = 2^{2}$				&\\
&$5 \Rightarrow 16 = 2^{4}$				&\\
&$21 \Rightarrow 64 = 2^{6}$			&\\
&$85 \Rightarrow 256 = 2^{8}$			&\\
&$341 \Rightarrow 1024 = 2^{10}$		&\\
\hline
&&\\
\multicolumn{3}{l}{All seeds are 4-form numbers, and successive seeds are determined by using} \\
\multicolumn{3}{l}{the previous seed as the index of a 4-form number.} \\

&$5 = 4(1) + 1$							&\\
&$21 = 4(5) + 1$							&\\
&$85 = 4(21) + 1$							&\\
&$341 = 4(85) + 1$						&\\
&$K_{i+1} = 4K_{i} + 1 $					& $	K \in N, K_{1} = 1$\\
\hline
\end{tabular}
\label{table:Seeds}
\end{table}

Cadogan identified seeds as \emph{Set A} in (2.6) of [2].

\section{Columns}
Another method of analyzing Collatz sequences is based on assigning the natural numbers into a table of 12 columns, with each column being the value of the number Mod 12 (column = C mod 12). We identify mod 12 results of zero as column 12 for ease of analysis and discussion.
 
The characteristics of a certain column of the table apply to all numbers in that column.

Applying the column concept to the standard number forms:
\newline\indent 2-form numbers (evens) fall in columns 2, 4, 6, 8, 10, and 12
\newline\indent4-form numbers (4n+1) fall in columns 1, 5, and 9
\newline\indent8- and higher form numbers (non-standard form 4n+3) fall in columns 3, 7, and 11.

We can see that the number forms are mutually exclusive and collectively exhaustive 
with respect to the 12 columns.

All natural numbers can be defined in \emph{column form} as $12r + column$, where r is the row of the table, starting with 0. For example, 27 is column form $12*2 + 3$.

Applying the appropriate Collatz step to each column form:

\begin{table}[!htbp]
\caption{Column Analysis of Collatz Steps} 			
\centering 										
\begin{tabular}{ccccc} 						
\hline\hline 									
\rule{0pt}{4mm}
			&Std or 			& Column-form	& Mixed Column Form	& Resulting \\
Column		&Non-Std Form	& Before Step		& After Step			& Column\\
\hline										
\rule{0pt}{4mm}

1			&	4n+1			&	12r+1			&	36r+4				&	4\\
2			&	2n 				& 	12r+2			&	6r+1				&	1 or 7\\
3			&	4n+3			&	12r+3 			&	36r+10				&	10\\
4			&	2n 				&	12r+4 			&	6r+2				&	2 or 8\\
5			&	4n+1			&	12r+5 			&	36r+16				&   4\\
6			&	2n 				&	12r+6 			&	6r+3				&	3 or 9\\
7			&	4n+3 			&	12r+7 			&	36r+22				&	10\\
8			&	2n 				&	12r+8 			&	6r+4				&	4 or 10\\
9			&	4n+1 			&	12r+9 			&	36r+28				&   4\\
10			&	2n 				&	12r+10 		&	6r+5				&	5 or 11\\
11			&	4n+3 			&	12r+11 		&	36r+34				&	10\\
12			&	2n 				&	12r+12 		&	6r+6				&	6 or 12\\

\hline 											
\end{tabular}
\label{table:ColumnAnalysis}
\end{table}

Note that the results of Collatz steps on all 4n+1 numbers end up in column 4 and the Collatz steps on all 4n+3 numbers end up in column 10.

\newpage
A map of the column analysis will help us understand the process. Numbers in the circles represent columns and arrows represent steps.

\vspace{-4mm}
\begin{figure}[!htbp]
\caption{Column Steps Map}
\begin{center}
\fbox{\includegraphics[width=7cm]{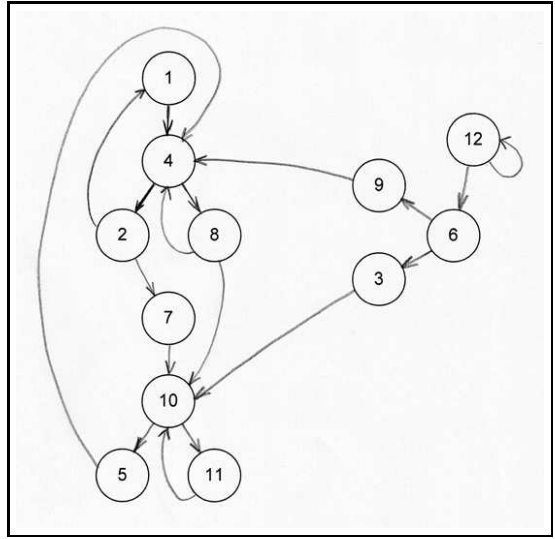}}
\end{center}
\label{figure:ColumnMap}
\end{figure}

\vspace{-4mm}
It’s clear from the diagram that multiples of 3 can only be initial values of Collatz sequences.

Applying the column concept to a 64-cascade:

\vspace{-4mm}
\begin{table}[!htbp]
\caption{Column Analysis of a General 64-Cascade} 			
\centering 										
\begin{tabular}{ccccc} 						
\hline\hline 									
\rule{0pt}{4mm}
&Composite & &&\\
Step&Form&Standard Form&Column Form&Column\\
\hline										
\rule{0pt}{4mm}

0	&	$64n+31$		&	$64n+31$			&	Indeterminate		&	3, 7, or 11		\\
1	&	$192n+94$		&	$2(96n+47)$		&$	12(16n+7)+10$	&$	10$				\\
2	&	$96n+47$		&	$32(3n+1)+15$	&$	12(8n+3)+11$		&$	11$				\\
3	&	$288n+142$	&	$2(144n+71)$		&$	12(24n+11)+10$	&$	10$				\\
4	&	$144n+71$		&	$16(9n+4)+7$		&$	12(12n+5)+11$	&$	11$				\\
5	&	$432n+214$	&	$2(216n+107)$	&$	12(36n+17)+10$	&$	10$				\\
6	&	$216n+107$	&	$8(27n+13)+3$	&$	12(18n+8)+11$	&$	11$				\\
7	&	$648n+322$	&	$2(324n+161)$	&$	12(54n+26)+10$	&$	10$				\\
8	&	$324n+161$	&	$4(81n+40)+1$	&$	12(27n+13)+5$	&$	5$				\\	
9	&	$972n+484$	&	$2(486n+242)$	&$	12(81n+40)+4$	&$	4$				\\
10	&	$486n+242$	&	$2(243n+121)$	&	indeterminate		&	2 or 8			\\
11	&	$243n+121$	&	indeterminate		&	indeterminate		&	1, 7, 4, or 10	\\

\hline 							
\end{tabular}
\label{table:64CascadeAnalysis}
\end{table}

We see that the bulk of the action of the cascade alternates between columns 10 and 11, and that the cascade ends in column 1, 7, 4, or 10 (the four 3t+1 columns). In fact, this is true for all 4n+3 cascades.

Applying the column concept to a plummet from a power of 2:

\begin{table}[!htbp]
\caption{Column Analysis of a Plummet} 			
\centering 										
\begin{tabular}{ccccc} 						
\hline\hline 									
\rule{0pt}{4mm}

Step	&Value					& Standard Form	&Column Form		&Column\\
\hline										
\rule{0pt}{4mm}

0		& $85 \text{(a seed)}$	& $4(21)+1$ 		& $	12(7)+1$ 		& $	1$\\
1		& $	256	$ 				& $2(128)$ 		& $	12(21)+4$ 		& $	4$\\
2		& $	128	$ 				& $2(64)$ 			& $	12(10)+8$ 		& $	8$\\
3		& $	64	$ 				& $2(32)$ 			& $	12(5)+4$ 		& $	4$\\
4		& $	32	$ 				& $2(16)$ 			& $	12(2)+8$ 		& $	8$\\
5		& $	16	$ 				& $2(8)$ 			& $	12(1)+4$ 		& $	4$\\
6		& $	8	$ 				& $2(4)$ 			& $	12(0)+8$ 		& $	8$\\	
7		& $	4	$ 				& $2(2)$ 			& $	12(0)+4$ 		& $	4$\\
8		& $	2	$ 				& $2(1)$ 			& $	12(0)+2$ 		& $	2$\\
9		& $	1	$ 				& $1	$ 			& $12(0)+1$ 		& $	1$\\

\hline 
\end{tabular}
\label{table:PlummetAnalysis}
\end{table}

In the case of a plummet, the bulk of the action takes place between columns 4 and 8, and the plummet always leads to a value of 1. The number of steps to reach a value of 1 is called the \emph{total stopping time}.

The maximum value in an odd cascade is the result of the odd step applied to the 4-form number near the end of the cascade (for example, after step 8 in Table \ref{table:64CascadeAnalysis}). The maximum value of the cascade is four times the value of the cascade transform and is always a column 4 number.

\section{Stopping Times}
\subsection{Stopping Times of Standard Forms}
Standard 2-form numbers have a stopping time of 1 because they immediately result in a smaller number.

\begin{table}[!htbp]
\caption{Stopping Times of Standard Forms}
\centering 										
\begin{tabular}{c l p{65mm}} 						
\hline\hline
\rule{0pt}{4mm}
Step	&   Form						&														\\
\hline
\rule{0pt}{4mm}
0		&	$2n$						&	Even 												\\
1		&	$2n/2=n$					&	Odd or even, depending on n, 					\\
     		&								&  but smaller value than $2n\text{ for all }n \in N$ 		\\

&& \\
\multicolumn{3}{l}{Standard 4-form numbers have a stopping time of 3.} \\

Step	&    Form						&														\\
0		&	$4n+1$		 			&	Odd 												\\
1		&	$3(4n+1)+1=12n+4$		&	Even												\\
2    	&	$(12n+4)/2=6n+2$		&   Even												\\			
3		&   $(6n+2)/2=3n+1$			&   Odd or even, depending on n,						\\
		&								&   but smaller value than $4n+1\text{ for all }n \in N$ 	\\
		&								&   ($n=0$ is the trivial case for the number 1) 	\\
\hline
\end{tabular}
\label{table:StdFormsStopTimes}
\end{table}

\newpage 

\begin{table}[!htbp]
\caption*{Table 14: Stopping Times of Standard Forms (continued)}
\centering 										
\begin{tabular}{c l p{5cm}} 						
\hline\hline
&& \\
\multicolumn{3}{p{12cm}}{Standard 8-form numbers are indeterminate after 5 steps (their parity depends on n), and as a group do not have a stopping time.}\\
&& \\Step& Form&\\
\hline
\rule{0pt}{4mm}
0	&	$8n+3$		 			&	Odd, 8-form 		\\
1	&	$3(8n+3)+1=24n+10$	&	Even		\\
2    &	$(24n+10)/2=12n+5$		& Odd, 4-form $=4(3n+1)+1$\\			
3	&   $3(12n+5)+1=36n+16$	& Even 	\\
4	&   $(36n+16)/2=18n+8$		&  Even 	\\
5	&   $(18n+8)/2=9n+4$		&  Odd or even, depending on n,	\\
	&   								&  but larger value than 		 	\\
	&								&  $8n+3\text{ for all }n \in W$ 	\\
&& \\
\multicolumn{3}{l}{Standard 32-form numbers are indeterminate after 9 steps (their parity}\\
\multicolumn{3}{l}{ depends on n), and as a group do not have a stopping time.}\\
								
Step	&    Form							&	\\
\hline
\rule{0pt}{4mm}
0	&	32n + 15							&	Odd, 32-form	\\
1	&	3(32n + 15) + 1 = 96n + 46		& Even	\\
2	&	(96n + 46) / 2 = 48n + 23			& Odd, 16-form $= 16(3n + 1) + 7$	\\
3	&	3(48n + 23) + 1 = 144n + 70		& Even	\\
4	&	(144n + 70) / 2 = 72n + 35		& Odd, 8-form $= 8(9n + 4) + 3$	\\
5	&	3(72n + 35) + 1 = 216n + 106	& Even	\\
6	&	(216n + 106) / 2 = 108n + 53		& Odd, 4-form $= 4(27n + 13) + 1$	\\
7	&	3(108n + 53) + 1 = 324n + 160	& Even	\\
8	&	(324n + 160) / 2 = 162n + 80		& Even	\\
9	&	(162n + 80) / 2 = 81n + 40		& Odd or even, depending on n,		\\ 
	&										& but larger value than 						\\
	&										& $32n + 15 \text{ for all }n \in W$	\\
\hline
&& \\
\multicolumn{3}{p{12cm}}{We find that 8- and higher standard forms, as a group, are indeterminate and do not have a stopping time because the base becomes odd before the stopping time of the offset is reached.} \\
&& \\
\multicolumn{3}{p{12cm}}{All 8- and higher standard forms are indeterminate after $2p-1$ steps (i.e, at the end of the cascade).} \\

\end{tabular}
\end{table}

\newpage
\subsection{Stopping Times of Composite Forms}
The odd/even parity of the result of each step in the Collatz sequence of a composite form is the same as the parity of the offset, as long as the base of the composite form at that step is even. When the base is odd (a mixed form), the parity of the number is indeterminate (depends on the value of n).

Standard form $8n+3$ is indeterminate after 5 steps, having a mixed form of $9n+4$, with a parity that depends on n and a value larger than $8n+3$ for all $n \in W$. The offset of this form, 3, has a stopping time of 6, which is found by applying the Collatz algorithm. If the base were even after step 5, the next step would result in an offset of 2, which is less than the starting offset. If we multiply the initial base (8) by 2, we get the composite form $16n+3$. Analyzing this form we find that it has a stopping time of 6 with a mixed form $9n+2$, which is smaller than $16n+3$ for all $n \in W$. 

Standard form $32n+15$ is indeterminate after 9 steps, having a mixed form of $81n+40$, with a parity that depends on n and a value larger than $32n+15$ for all $n \in W$. The offset of this form, 15, has a stopping time of 11, which is found by applying the Collatz algorithm. If the base were even through the step just prior to the stopping time of the offset (in this case step 10), the next step would result in an offset of 10, which is less than the starting offset. If we multiply the initial base (32) by 4 (we’ll describe how to determine this factor later), we get the composite form $128n+15$. Analyzing this form (see Table \ref{table:128CascadeAnalysis}) we find that it has a stopping time of 11 with a mixed form $81n+10$, which is smaller than $128n+15 \text{ for all }n \in W$. 

\begin{table}[!htbp]
\caption{Collatz Sequence for Composite Form $128n+15$} 			
\centering 										
\begin{tabular}{c l p{5cm}} 						
\hline\hline
\rule{0pt}{4mm}

Step&	Form	\\
0	&	$128n + 15$							& Odd, 32-form $= 32(4n) + 15$	\\
1	&	$3(128n + 15) + 1 = 384n + 46$	& Even $(\#1)$	\\
2	&	$(384n + 46) / 2 = 192n + 23$		& Odd, 16-form $= 16(12n + 1) + 7$	\\
3	&	$3(192n + 23) + 1 = 576n + 70$	& Even $(\#2)$	\\
4	&	$(576n + 70) / 2 = 288n + 35$		& Odd, 8-form $= 8(36n + 4) + 3$	\\
5	&	$3(288n + 35) + 1 = 864n + 106$	& Even $(\#3)$	\\
6	&	$(864n + 106) / 2 = 432n + 53$		& Odd, 4-form $= 4(108n + 13) + 1$	\\
7	&	$3(432n + 53) + 1 = 1296n + 160$	& Even $(\#4)$	\\
8	&	$(1296n + 160) / 2 = 648n + 80$		& Even $(\#5)$	\\
9	&	$(648n + 80) / 2 = 324n + 40$		& Even $(\#6)$\\
10	&	$(324n + 40) / 2 = 162n + 20$		& Even $(\#7)$	\\
11	&	$(162n + 20) / 2 = 81n + 10$			& Odd or even, depending on n, 	\\
	&											& but smaller value than \\
	&											& $128n + 15 \text{ for all } n \in W$	 \\
\hline
\end{tabular}
\label{table:128CascadeAnalysis}
\end{table}
\vspace{-2mm}
\textbf{In order for a number in a composite form sequence to be smaller than the starting number, the power of 2 in the base of the starting number must be at least equal to the number of even steps in the Collatz sequence of the offset through its stopping time.}

Looking at the sequence for composite form $128n+15$ (Table \ref{table:128CascadeAnalysis}), the number of even steps until the resulting offset is less than the starting offset is 7. This means that the smallest base for a composite form with an offset of 15 that ensures a finite stopping time is $2^{7} = 128$. All composite forms with an offset of 15 and a power of 2 greater than 128 in the base also have a stopping time of 11. 

This is true for all sequences. As long as the power of 2 in the base is equal to or larger than the number of even steps in the Collatz sequence of the offset through its stopping time, the sequence has a finite stopping time. \textbf{This means that every number represented by a composite form that has a stopping time also has that stopping time.} This was also identified by Garner in [3].

\subsection{Calculating the Minimum Base That Ensures a Stopping Time for a Particular Offset}
There is a direct relationship between stopping time and number of even steps to the stopping time in a Collatz sequence. We can, in fact, calculate the minimum power of 2 in the base that is required to ensure that a composite form will have a finite stopping time, based on the stopping time of the offset.

In order to reach a number smaller than the starting number, the ratio of the effect of odd steps (multiply by 3) to even steps (divide by 2) must be between 0.5 and 1. (We can disregard the $+ 1$ term in the odd steps for now – it doesn’t seem to affect these calculations.)
This means that

\begin{table}[!htbp]
\caption{Calculating Minimum Base to Ensure Stopping Time}
\centering 										
\begin{tabular}{l p{5cm}r} 						
\hline\hline
\rule{0pt}{4mm}
$0.5 < 3^{S-E}/2^{E} < 1$				& where S is the stopping time of the offset, E is the number of even steps to the stopping time, and $S-E$ is the number of odd steps to the stopping time.&(5)\\
$0.5 < 3^{S}/(3^{E}*2^{E}) < 1$	& rearrange $3^{-E}$ &\\
$0.5 < 3^{S}/6^{E} < 1$				& combine $3^{E} *2^{E}$&\\
$2 > 6^{E}/3^{S} > 1$				& take reciprocals&\\
$1 < 6^{E}/3^{S} < 2$				& reorder&\\
$3^{S} < 6^{E} < 2*3^{S}$			& multiply by $3^{S}$&\\
$S*ln(3) < E*ln(6) < ln(2) + S*ln(3)$	& take natural logs&\\
										&										&\\
\multicolumn{2}{p{10cm}}{$S*ln(3)/ln(6) < E < S*ln(3)/ln(6)+ln(2)/ln(6)$}&(6)\\
										& divide by $ln(6) \approx 1.792$		&\\
\multicolumn{2}{p{10cm}}{or approximately $0.613*S < E < 0.613*S+0.387$}	&\\
\end{tabular}
\label{table:CalcMinimumBase}
\end{table}

Since E must be a natural number there is at most one value of E that will satisfy inequality (6) for any value of S; the integer portion of $0.613*S$ must be smaller than the integer portion of $0.613*S + 0.387$. \newpage We can calculate E for a given value of S using the following equation:

\begin{table}[!htbp]
\centering 										
\begin{tabular}{p{10cm}r} 						

$E=[int(S*A+B)-int(S*A)]*int(S*A+B)$    &(7)\\
where $A=ln(3)/ln(6)$ and $B=ln(2)/ln(6)$&\\

\end{tabular}
\end{table}

For many values of S there is no value of E that satisfies inequality (6). This should mean that there is no sequence that has that value of S as its stopping time, and the authors have found this to be true for all numbers tested (for example, when S = 4, 5, or 7).

Table \ref{table:EvensToStopTime} shows number of even steps to the stopping time and power of 2 in the base to ensure composite forms have stopping times for some small values of S. Numbers omitted from the stopping time column do not exist as the stopping time of any Collatz sequence. (Identified by Terras with different nomenclature, Table B, p. 248 [6].)

\begin{table}[!htbp]
\caption{Even Steps Required to Stopping Time} 			
\centering 										
\begin{tabular}{ccc} 						
\hline\hline
\rule{0pt}{4mm}
				&							&Minimum\\
Stopping Time	&	Even Steps Required	&	Value of Base	\\
S				&	E						&	$2^{E}$	\\
\hline
\rule{0pt}{4mm}
  1				&	  1						&	2	\\
  3				&	  2						& 	4	\\
  6				& 	  4						&	16	\\
  8				&	  5						&	32	\\
11				&	  7						&	128	\\
13				&	  8						&	256	\\
16				&	10						&	1024	\\
19				&	12						&	4096	\\
21				&	13						&	18192	\\
24				&	15						&	32768	\\
26				&	16						&	65536	\\
…				&	…						&	…	\\

\end{tabular}
\label{table:EvensToStopTime}
\end{table}

For example, this means that all sequences that have a stopping time of 8 are part of a group of numbers that have a base of 32. It turns out that only two composite forms have a stopping time of 8: $32n+11$ and $32n+23$.

This also allows us to determine the group in which a given offset will have a finite stopping time. For example, the Collatz sequence for 27 has a stopping time of 96. When we calculate E for $S=96$ using equation (7), we get a value of 59, which is the number of even steps in the sequence for 27 until the stopping time.  Thus, all numbers of the composite form $2^{59}n+27$ have a stopping time of 96.

\newpage
\subsection{Principal Composite Forms}
As a group, no standard form higher than 4 has a stopping time.

As we analyze larger and larger numbers, we find that most numbers can be described by composite forms having finite stopping times, and offsets that are smaller than the number being analyzed. For example, 41 is part of group $4n+1$ with a stopping time of 3, and 43 is part of group $32n+11$ with a stopping time of 8 [3]. (Kannan and Moorthy [4] identified similar results, with some of the results in other than lowest form.) 

However, some numbers are not part of a group with a stopping time, and an offset smaller than the number, such as 47, which is part of group $2^{54}n+47$; or 27, which is part of group $2^{59}n+27$. We call a composite form that has a stopping time, and an offset equal to the number being analyzed a \emph{principal form}, because no numbers smaller than the offset are described by that form.

Table 12 shows principal forms required to describe all numbers up to 100.

\begin{table}[!htbp]
\caption{Principal Forms} 			
\centering 										
\begin{tabular}{c} 						
 required to describe \\
 all numbers up to 100\\
\hline
\rule{0pt}{4mm}

$2 * n$	\\
$4 * n + 1$	\\
$16 * n + 3$	\\
$128 * n + 7$	\\
$32 * n + 11$	\\
$128 * n + 15$	\\
$32 * n + 23$	\\
$2^{59} * n + 27$	\\
$2^{56} * n + 31$	\\
$256 * n + 39$	\\
$2^{54} * n + 47$	\\
$128 * n + 59$	\\
$2^{54} * n + 63$	\\
$2^{51} * n + 71$	\\
$256 * n + 79$	\\
$2^{45} * n + 91$	\\
$256 * n + 95$	\\
\end{tabular}
\label{table:PrincipalForms}
\end{table}

An additional 10 principal forms are required to describe all numbers between 101 and 200, and 9 more are required for all numbers between 201 and 300. Garner identified a number of these Principal Composite Forms in [3].

If there exists a number M such that no additional principal forms are required to ensure a stopping time for all numbers greater than M, the Collatz Conjecture can be proven true. We find, however, that the number of principal forms per 10,000 numbers analyzed behaves roughly as shown in Figure 2. 

\begin{figure}
\caption{Average Principal Forms per 10,000 Numbers}
\begin{center}
\fbox{\includegraphics[width=10cm]{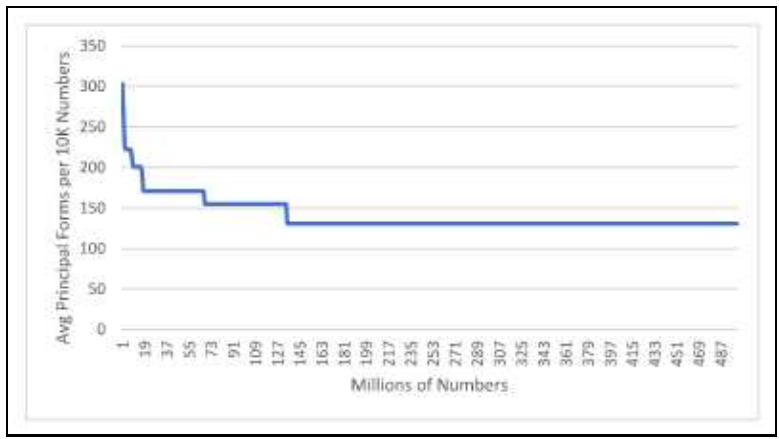}}
\end{center}
\label{figure:FormsPer10K}
\end{figure}

\newpage
The line in Figure \ref{figure:FormsPer10K} appears to reach a constant value, but there is actually some variation in the average values of principal forms per 10,000 numbers tested. Table 19 shows statistics for some interesting points in Figure \ref{figure:FormsPer10K}.

\begin{table}[!htbp]
\caption{Principal Forms per 10,000 Numbers} 			
\centering 										
\begin{tabular}{ccccc} 						
\hline\hline
\rule{0pt}{4mm}
Million*	&Average	&Std. Dev.	&	Max	&	Min\\
\hline
\rule{0pt}{4mm}

  1		&	303.4	&	43.6	&	589		&	268	\\
  2		&	261.7	&	9.2		&	297		&	245	\\
  3		&	225.9	&	14.5	&	270		&	200	\\
  4		&	222.4	&	7.4		&	239		&	197	\\
  5		&	222.2	&	8.0		&	246		&	196	\\
  6		&	222.2	&	7.8		&	237		&	199	\\
  7		&	222.4	&	7.3		&	240		&	201	\\
  8		&	222.2	&	8.1		&	242		&	201	\\
  9		&	209.2	&	13.1	&	236		&	183	\\
10		&	201.4	&	8.0		&	216		&	182	\\
…		&			&			&			&		\\	
16		&	201.3	&	7.0		&	216		&	186	\\
17		&	194.8	&	15.3	&	216		&	153	\\
18		&	170.6	&	8.6		&	195		&	148	\\
19		&	170.5	&	7.4		&	186		&	148	\\
…		&			&			&			&		\\			
66		&	170.4	&	6.8		&	186		&	152	\\
67		&	170.7	&	7.0		&	187		&	151	\\
68		&	157.0	&	9.8		&	190		&	135	\\
69		&	153.6	&	7.6		&	170		&	135	\\
70		&	154.4	&	9.2		&	177		&	133	\\
…		&			&			&			&		\\				
133		&	154.6	&	7.9		&	178		&	135	\\
134		&	154.1	&	7.6		&	174		&	135	\\
135		&	137.3	&	12.8	&	172		&	113	\\
136		&	129.3	&	8.3		&	151		&	112	\\
…		&			&			&			&		\\				
250		&	131.5	&	7.9		&	149		&	115	\\
…		&			&			&			&		\\				
500	&		131.4	&	8.2		&	150		&	110	\\
\multicolumn{5}{l}{* the one million numbers ending at}\\
\multicolumn{5}{l}{this number times 1,000,000}\\

\end{tabular}
\label{table:FormsPer10K}
\end{table}

This appears to be a downward trend which may continue, indicating that the number M may exist, or may reach a limit of some number of principal forms per 10,000 numbers which would indicate that this approach cannot be used to confirm the Collatz Conjecture. 

\indent
Continuing analysis of Collatz sequences of selected groups of very large numbers shows that principal composite forms decrease to 0 or 1 per 10,000 for those groups of numbers. The authors have found an astonishing string of 50,000 consecutive odd numbers with the same total stopping time (number of steps to reach 1). The numbers are in the range $10^{142} - 10^{6} + 1$ through 99,999. Analysis of principal composite forms for very large numbers continues.

\subsection{Generating and Analyzing Composite Forms}
We can generate composite forms by replacing the index of a standard form with expressions for $\#$-forms. For example, if we replace the index n in the standard 8-form with the 2-form expression ($2n$) we get a new expression for an 8-form with a 2-form base: $8(2n) + 3 = 16n + 3$. We identify this composite form as $8.2$. Another example: form $16.4.8$ is $16(4(8n + 3) + 1) + 7 = 512n + 215$.

As stated earlier in this paper, the form of the number that results from a cascade determines whether the sequence will go up or down and by how much. If we could determine the forms of cascade results we might gain some insight into the overall behavior of Collatz sequences.

The 4-cascade transforms the starting value of $4n+1$ into a cascade result of $3n+1$. The values of $3n+1$ are a \emph{mix} of forms that depend on the value of n. The 4.2-cascade transform is $8n+1$ to $6n+1$, which is also a mix of different forms, depending on the value of n.

However, the 4.4 and higher 4.x cascade transforms are all 2-form numbers. Form 4.4 can be described as $4(4n+1)+1$, and a cascade transforms it to $3(4n+1)+1$, which is $12n+4$, which is always even (2-form). Form 4.8 and higher can be described as $4(4n+3)+1$, and a cascade transforms it to $3(4n+3)+1$, which is $12n+10$, which is always even (2-form).

Tables \ref{table:Level1} through \ref{table:Level3} show forms of cascade results for a number of composite forms. Note that each row of the tables is one of two general types:

\begin{adjustwidth}{1cm}{0cm}

\hspace{5 mm}For Type 1, the value shown in the .2 column is \emph{Mix}, followed by a particularform for .4 and higher columns. 
 
For Type 2, the value in the .2 column is a particular form, followed by increasing forms until a \emph{Mix} results, followed by the next higher form after the one just preceding the \emph{Mix} for all higher columns.

\end{adjustwidth}
In each row of each table, the smallest form of numbers in the mix column is 1 more than the highest form in the other columns of the row. The forms for the \emph{Mix} values are shown in the rightmost column of each row. 

Tables \ref{table:Level2} and \ref{table:Level3} show results for the \emph{Mix} values from the previous table.

While these results are not conclusive, they indicate yet another interesting pattern in Collatz sequences. 

\begin{table}[!htbp]
\caption{Standard Forms of Cascade Results - Level 1} 			
\centering 										
\begin{tabular}{ccccccccc} 						
\multicolumn{9}{c}{(Standard forms of numbers resulting from applicable cascades)}\\
\hline\hline
\rule{0pt}{4mm}
&\multicolumn{8}{c}{<========== Level 2 ==========>}\\
Level 1&	.2&	.4&	.8&	.16&	.32&	.64&	.128&	Mix\\
\hline
\rule{0pt}{4mm}
4			&Mix	&2		&2		&2		&2		&2		&2		&4+  \\
8			&2		&4		&Mix	&8		&8		&8		&8		&16+  \\
16			&Mix	&2		&2		&2		&2		&2		&2		&4+  \\
32			&2		&4		&8		&Mix	&16	&16	&16	&32+  \\
64			&Mix	&2		&2		&2		&2		&2		&2		&4+  \\	
128			&2		&4		&Mix	&8		&8		&8		&8		&16+  \\
256			&Mix	&2		&2		&2		&2		&2		&2		&4+  \\
512			&2		&4		&8		&16	&Mix	&32	&32	&64+  \\
1024		&Mix	&2		&2		&2		&2		&2		&2		&4+  \\
2048		&2		&4		&Mix	&8		&8		&8		&8		&16+  \\
4096		&Mix	&2		&2		&2		&2		&2		&2		&4+  \\
8192		&2		&4		&8		&Mix	&16	&16	&16	&32+  \\
16384		&Mix	&2		&2		&2		&2		&2		&2		&4+  \\
32768		&2		&4		&Mix	&8		&8		&8		&8		&16+  \\
65536		&Mix	&2		&2		&2		&2		&2		&2		&4+  \\	
131072		&2		&4		&8		&16	&32	&Mix	&64	&128+  \\
262144		&Mix	&2		&2		&2		&2		&2		&2		&4+  \\
524288		&2		&4		&Mix	&8		&8		&8		&8		&16+  \\
1048576	&Mix	&2		&2		&2		&2		&2		&2		&4+  \\
2097152	&2		&4		&8		&Mix	&16	&16	&16	&32+  \\
\end{tabular}
\label{table:Level1}
\end{table}

\begin{table}[!htbp]
\caption{Standard Forms of Cascade Results - Level 2} 			
\centering 										
\begin{tabular}{ccccccccc} 						
\multicolumn{9}{c}{(Standard forms of numbers resulting from applicable cascades)}\\
\multicolumn{9}{c}{Entries shown only for \emph{Mix} Values in Table \ref{table:Level1}}\\
\hline\hline
\rule{0pt}{4mm}
&\multicolumn{8}{c}{<========== Level 3 ==========>}\\
Level 2			&	.2	&	.4	&	.8	&	.16	&	.32	&	.64	&	.128&	Mix\\
\hline
\rule{0pt}{4mm}
4.2				&4		&Mix	&8		&8		&8		&8		&8		&16+ \\
8.8				&Mix	&16	&16	&16	&16	&16	&16	&32+ \\
16.2			&4		&8		&Mix	&16	&16	&16	&16	&32+ \\
32.16			&Mix	&32	&32	&32	&32	&32	&32	&64+ \\
64.2			&4		&Mix	&8		&8		&8		&8		&8		&16+ \\
128.8			&16	&Mix	&32	&32	&32	&32	&32	&64+ \\
256.2			&4		&8		&16	&Mix	&32	&32	&32	&64+ \\
512.32			&Mix	&64	&64	&64	&64	&64	&64	&128+ \\
1024.2			&4		&Mix	&8		&8		&8		&8		&8		&16+ \\
2048.8			&Mix	&16	&16	&16	&16	&16	&16	&32+ \\
4096.2			&4		&8		&Mix	&16	&16	&16	&16	&32+ \\
8192.16		&32	&64	&128	&Mix	&256	&256	&256	&512+ \\
16384.2		&4		&Mix	&8		&8		&8		&8		&8		&16+ \\
32768.8		&16	&32	&Mix	&64	&64	&64	&64	&128+ \\
65536.2		&4		&8		&16	&32	&Mix	&64	&64	&128+ \\
131072.64		&Mix	&128	&128	&128	&128	&128	&128	&256+ \\
262144.2		&4		&Mix	&8		&8		&8		&8		&8		&16+ \\
524288.8		&Mix	&16	&16	&16	&16	&16	&16	&32+ \\
1048576.2		&4		&8		&Mix	&16	&16	&16	&16	&32+ \\
2097152.16	&Mix	&32	&32	&32	&32	&32	&32	&64+ \\
\end{tabular}
\label{table:Level2}
\end{table}

\begin{table}[!htbp]
\caption{Standard Forms of Cascade Results - Level 3} 			
\centering 										
\begin{tabular}{ccccccccc} 						
\multicolumn{9}{c}{(Standard forms of numbers resulting from applicable cascades)}\\
\multicolumn{9}{c}{Entries shown only for \emph{Mix} Values in Table \ref{table:Level2}}\\
\hline\hline
\rule{0pt}{4mm}
&\multicolumn{8}{c}{<========== Level 4 ==========>}\\
Level 3			&.2		&.4		&.8		&.16	&.32	&.64	&.128	&Mix\\
\hline
\rule{0pt}{4mm}
4.2.4			&16	&Mix	&32	&32	&32	&32	&32	&64+	\\
8.8.2			&Mix	&32	&32	&32	&32	&32	&32	&64+	\\
16.2.8			&32	&64	&128	&256	&Mix	&512	&512	&1024+	\\
32.16.2			&64	&Mix	&128	&128	&128	&128	&128	&256+	\\
64.2.4			&Mix	&16	&16	&16	&16	&16	&16	&32+	\\
128.8.4			&Mix	&64	&64	&64	&64	&64	&64	&128+	\\
256.2.16		&64	&128	&Mix	&256	&256	&256	&256	&512+	\\
512.32.2		&128	&256	&Mix	&512	&512	&512	&512	&1024+	\\
1024.2.4		&16	&32	&Mix	&64	&64	&64	&64	&128+	\\
2048.8.2		&32	&64	&Mix	&128	&128	&128	&128	&256+	\\
4096.2.8		&Mix	&32	&32	&32	&32	&32	&32	&64+	\\
8192.16.16		&512	&Mix	&1024	&1024	&1024	&1024	&1024	&2048+	\\
16384.2.4		&Mix	&16	&16	&16	&16	&16	&16	&32+	\\
32768.8.8		&Mix	&128	&128	&128	&128	&128	&128	&256+	\\
65536.2.32		&128	&256	&Mix	&512	&512	&512	&512	&1024+	\\
131072.64.2	&256	&512	&1024	&Mix	&2048	&2048	&2048	&4096+	\\
262144.2.4		&16	&Mix	&32	&32	&32	&32	&32	&64+	\\
524288.8.2		&Mix	&32	&32	&32	&32	&32	&32	&64+	\\
1048576.2.8	&32	&Mix	&64	&64	&64	&64	&64	&128+	\\
2097152.16.2	&Mix	&64	&64	&64	&64	&64	&64	&128+	\\
&&&&&&&&\\
\multicolumn{9}{c}{  It appears that this pattern continues at higher levels of detail.}\\
\end{tabular}
\label{table:Level3}
\end{table}

\newpage
\section{Conclusions}
Forms help us analyze large groups of numbers having the same characteristics, and cascades help us see the rigid structure within Collatz sequences, providing important insights into what appear to be random ups and downs. Cascades also help us understand the next steps after a cascade – namely the start of another cascade. Understanding results of cascades is key to understanding the Collatz Conjecture. Unfortunately, any cascade can result in a number of any form. 

Column analysis exposes hidden patterns within cascades.

Successive reverse cascades of non-multiples of 3 all seem to lead to an odd multiple of 3, providing another interesting approach to possibly confirming the Collatz Conjecture.

The relationship between stopping time and number of even steps to stopping time allows us to group many numbers having the same principal composite form – all of which have the same stopping time. Analysis of principal composite forms for very large numbers continues.

Total stopping time is not related to the form of a number or its stopping time.

The authors hope that the patterns we have discovered in Collatz sequences may lead to a proof of the Collatz Conjecture. 

\section{References}
[1] Barina, David (2020). "Convergence verification of the Collatz problem." The Journal of Supercomputing.  doi:10.1007/s11227-020-03368-x. S2CID 220294340.\vspace*{3mm}\newline
[2] Cadogan, Charles C (1991). "Some Observations on the 3x + 1 Problem.” Proc. Sixth Caribbean Conference on Combinatorics $\&$ Computing, University of the West Indies: St. Augustine Trinidad (C. C. Cadogan, Ed.) Jan 1991, 84-91.\vspace*{3mm}\newline
[3] Garner, Lynn E. (1981). “On the Collatz 3n+1 Algorithm.” Proceedings of the American Mathematical Society, Volume 82, Number 1, May 1981, 19-22.\vspace*{3mm}\newline
[4] Kannan, T. and Moorthy, C. Ganesa (2016). “Two Methods to Approach Collatz Conjecture.” International Journal of Mathematics And its Applications, Volume 4, Issue 1-B (2016), 161-167.\vspace*{3mm}\newline
[5] Mehendale, Dhananjay P. (2005). “Some Observations on the 3x+1 Problem.” arXiv:math/0504355 [math.GM] 18 Apr 2005.\vspace*{3mm}\newline
[6] Terras, Riho (1976). “A stopping time problem on the positive integers.” Acta Arithmetica 30 (1976), 241-252.\vspace*{3mm}\newline
[7] Trümper, Manfred (2006). “Handles, Hooks, and Scenarios: A fresh Look at the Collatz Conjecture.” arXiv:Math/0612228v1 [math.GM] 9 Dec 2006.

\newpage
\section{Appendices}
\subsection*{Appendix I - Proof of Lemma 1}
\textbf{Lemma 1:} All natural numbers can be expressed in standard form by the equation:

\begin{table}[!htbp]
\centering 		
\begin{tabular}{l l r}
$C=2^{p} n+2^{p-1}-1;$ 		& $p \in N, n \in W (N \cup \{0\})$  		&  (1)\\ 
\end{tabular}
\end{table}
This equation was developed independently by the authors, and was previously identified by Mehendale as Obs. 2 in [5] and by Cadogan as (3.2) of [2], and also describes even numbers when p=1.

\vspace{4pt}
\textbf{Proof of Lemma 1:}

\begin{table}[!htbp]
\begin{tabular}{lllr} 						
\multicolumn{4}{l}{All natural numbers are either even or odd.}\\
\multicolumn{4}{l}{Case 1: Even natural numbers $C = 2t$; $t \in N$,}\\
&	$C = 2t$						&	Standard 2-form,			&	 (8)	\\	
&	$C = 2t + 1 - 1$				&								&	\\
&	$C = 2^{p}n + 2^{p-1} - 1$ 	&	with $p = 1$ and $n = t$	&		\\	
			
\multicolumn{4}{l}{All even natural numbers are standard form $2^{p}n + 2^{p-1} - 1$}\\
\multicolumn{4}{l}{with $p = 1$ and $n \in N$.}\\
&&&\\

\multicolumn{4}{l}{Case 2: Odd natural numbers $C = 2t+1$; $t \in W$,}\\
&	$C = 2t+1$					&	Non-standard 2-form,	&	(9)	\\	
&	$C = 2t + 2 - 1$				&							&	\\
&	$C = 2^{p}n + 2^{p} - 1$ 	&	with $p = 1$ and $n = t$	&		\\	
			
\multicolumn{4}{l}{All natural numbers are either standard form ($2n$, even) }\\
\multicolumn{4}{l}{or non-standard 2-form ($2n+1$, odd).}\\
&&&\\

\multicolumn{4}{l}{Continuing, $t$ in equation (9) must be even or odd}	\\
\multicolumn{2}{l}{If $t$ is even, let $t=2m$}	&$m \in W$ 		&		\\
From (9)&	$C = 2(2m)+1$		&									&		\\	
&	$C = 4m+1$					&	Standard 4-form				& (10)	\\	
&	$C = 4m+2-1$					&									&		\\
&	$C = 2^{p}n + 2^{p-1} - 1$ 	&	with $p = 2$ and $n = m$	&		\\	
&&&\\			
\multicolumn{2}{l}{If $t$ is odd, let $t=2m+1$}	&$m \in W$ 		&		\\
From (9)&	$C = 2(2m+1)+1$		&									&		\\	
&	$C = 4m+2+1$				&									&		\\	
&	$C = 4m+3$					&	Non-standard 4-form			& (11)	\\
&	$C = 4m+4-1$					&									&	\\	
&	$C = 2^{p}n + 2^{p} - 1$ 	&	with $p = 2$ and $n = m$	&		\\	
		
\multicolumn{4}{l}{All natural numbers of non-standard 2-form ($2n+1$, odd) are either}\\
\multicolumn{4}{l}{standard 4-form ($4n+1$) or non-standard 4-form ($4n+3$).}\\
&&&\\	
\end{tabular}
\end{table}

\begin{table}[!htbp]
\centering 										
\begin{tabular}{lllr} 						
\multicolumn{4}{l}{Continuing, $m$ in equation (11) must be even or odd}	\\
\multicolumn{2}{l}{If $m$ is even, let $m=2q$}	&$q \in W$ 		&		\\
From (11)&	$C = 4(2q)+3$		&									&		\\	
&	$C = 8q+3$					&	Standard 8-form				& (12)	\\	
&	$C = 8q+4-1$					&									&		\\
&	$C = 2^{p}n + 2^{p-1} - 1$ 	&	with $p = 3$ and $n = q$	&		\\	
&&&\\			
\multicolumn{2}{l}{If $m$ is odd, let $m=2q+1$}&$q \in W$ 		&		\\
From (11)&	$C = 4(2q+1)+3$		&									&		\\	
&	$C = 8q+4+3$					&									&		\\	
&	$C = 8q+7$					&	Non-standard 8-form			& (13)	\\
&	$C = 8q+8-1$					&									&	\\	
&	$C = 2^{p}n + 2^{p} - 1$ 	&	with $p = 3$ and $n = q$	&		\\	
			
\multicolumn{4}{l}{All natural numbers of non-standard 4-form ($4n+3$) are either}\\
\multicolumn{4}{l}{standard 8-form ($8n+3$) or non-standard 8-form ($8n+7$).}\\
&&&\\

\multicolumn{4}{l}{Continuing, $q$ in equation (13) must be even or odd}	\\
\multicolumn{2}{l}{If $q$ is even, let $q=2r$}	&$r \in W$ 		&		\\
From (13)&	$C = 8(2r)+7$		&									&		\\	
&	$C = 16r+7$					&	Standard 16-form				& (14)	\\	
&	$C = 16r+8-1$					&									&		\\
&	$C = 2^{p}n + 2^{p-1} - 1$ 	&	with $p = 4$ and $n = r$	&		\\	
&&&\\			
\multicolumn{2}{l}{If $q$ is odd, let $q=2r+1$}&$r \in W$ 		&		\\
From (13)&	$C = 8(2r+1)+7$		&									&		\\	
&	$C = 16r+8+7$					&									&		\\	
&	$C = 16r+15$					&	Non-standard 16-form			& (15)	\\
&	$C = 16r+16-1$					&									&	\\	
&	$C = 2^{p}n + 2^{p} - 1$ 	&	with $p = 4$ and $n = r$	&		\\	
		
\multicolumn{4}{l}{All natural numbers of non-standard 8-form ($8n+7$) are either}\\
\multicolumn{4}{l}{standard 16-form ($16n+7$)or non-standard 16-form ($16n+15$).}\\
&&&\\	
\end{tabular}
\end{table}

\newpage
In general, any number of a non-standard form $2^{p}n + 2^{p} - 1$ is either the next-higher standard form $2^{p+1}n + 2^{p} - 1$ or the next higher non-standard form $2^{p+1}n + 2^{p+1} - 1$.
 
This process continues indefinitely (p increases with each step) and will generate all odd numbers in the form $2^{p}n + 2^{p-1} - 1$ (standard form).

Standard 2-form and non-standard 2-form include all natural numbers (evens and odds, respectively). Similarly, standard 4-form and non-standard 4-form include all odd numbers (developed from equation 9). Non-standard 4-form includes all odd numbers that are not standard 4-form (numbers of a higher form), and non-standard 8-form includes all odd numbers that are not standard 4-form or standard 8-form. This continues to infinity, as illustrated below.\vspace{12pt}
\newline\noindent
|-------------- even form ------------------|--------------------- odd forms ------------------|\newline
|------------------- 2n ----------------------|------------------------ 2n+1 ---------------------|\newline
.\hspace{54mm}|----------- 4n+1 ----------|-------- 4n+3 --------|\newline
.\hspace{90mm}|-- 8n+3 -|-- 8n+7-|  \newline\vspace{12pt}
…and so forth, to infinity.

All odd numbers can be described in the form $2^{p}n + 2^{p-1} - 1$, with $p \geq 2$ and $n \in W$.

\vspace{18pt}
\textbf{Thus, Lemma 1 is proved: All natural numbers can be described in the form $2^{p}n + 2^{p-1} - 1$, with $p \in N$,  $p \geq 1$ and $n \in W$.}

\subsection*{Appendix II - Determining the Standard Form of a Number}

The standard form of a particular number can be calculated as follows:
\begin{table}[!htbp]
\centering 										
\begin{tabular}{lllr} 						
\multicolumn{4}{l}{Starting with equation (1), repeated here for reference:}			\\
&\multicolumn{2}{l}{$C=2^{p}n + 2^{p-1} - 1$}								& (1)	\\
&	$C+1 = 2^{p}n + 2^{p-1}$ 			&	Add 1								&		\\	
&	$C+1 = 2^{p-1}(2n+1)$				&	Factor out $2^{p-1}$				& (16)	\\	
&	$2^{p-1}=gcd(C+1,2^{30})$			&	gcd is greatest common divisor;	& 		\\
&											&	$2^{30}$ is simply a				&		\\
&											&	sufficiently large power of 2		&		\\
&	$p-1 = log_{2}(gcd(C+1,2^{30}))$ 		&	Take $log_{2}$ of each side		&		\\	
&	$p = log_{2}(gcd(C+1,2^{30}))+1$ 	&	Add 1								&(17)	\\	
From (16)&$2n+1=(C+1)/2^{p-1}$ 		& also identified by Cadogan [2] 		& 		\\
&	$2n=(C+1)/2^{p-1}-1$ 				& Subtract 1 							& 		\\
&	$2n=(C+1-2^{p-1})/2^{p-1}$		& Rearrange 							& 		\\
&	$n=(C+1-2^{p-1})/2^{p}	$			& Divide by 2 							& (18) \\
\end{tabular}
\end{table}

\newpage
\subsection*{Appendix III - Proof of Lemma 2}

\textbf{Lemma 2:} Every natural number can be described in standard form by a unique combination of p and n. Conversely, the combination of p and n for each natural number is unique.\vspace{18pt}

\textbf{Proof of Lemma 2:}\newline
Assume that there is another pair of numbers, say q and m, that generates the same natural number as p and n, respectively.

\begin{table}[!htbp]
\centering 										
\begin{tabular}{lllr} 						
Then,		&   $2^{p}n + 2^{p-1} - 1 = 2^{q}m + 2^{q-1} - 1$	&Assume $q \neq p$ and $m \neq n$			&\\
			&	$2^{p}n + 2^{p-1} = 2^{q}m + 2^{q-1}$			&												&\\
			&	$2^{p-1}(2n+1) = 2^{q-1}(2m+1)$					&	Factor out $2^{p-1}$ and $2^{q-1}$		&\\
			&	$2^{p-q} = (2m+1)/(2n+1)$							&												& (19)\\
But, 		&\multicolumn{3}{l}{$(2m+1)/(2n+1)$ must be odd, since both $(2m+1)$ and}\\
			&\multicolumn{3}{l}{$(2n+1)$ are odd}\\
Therefore,	&$2^{p-q}$ must also be odd								&												&	\\
&\multicolumn{3}{l}{The only power of 2 that is odd is $2^{0}=1$, so $p-q$ must equal 0,}\\
&\multicolumn{3}{l}{and q must equal p}\\
Thus,		&$2^{0}=1=(2m+1)/(2n+1)$ 							& 												& \\
from (19)	&															&												&\\
			&$2n+1=2m+1$											&												&\\
And		&$n=m$													&												&\\
&&&\\
\multicolumn{4}{l}{Since q must equal p and m must equal n, the original assumption is false.}\\
\end{tabular}
\end{table}

\textbf{Thus, Lemma 2 is proved: Each unique combination of p and n generates a unique natural number, and vice versa.}

\subsection*{Appendix IV - Example of an Odd Cycle}
\begin{table}[!htbp]
\centering 										
\begin{tabular}{llll} 						
&	\multicolumn{3}{l}{Applying an odd cycle to 87:}						\\
$C_{i}$			&   $ 87$				& Odd, standard 16-form	& $=16(5)+7$		\\
$C_{i+1}$		&   $ 3*87+1=262$	& Even, standard 2-form	& $=2(131)$		\\
$C_{i+2}$		&   $ 262/2=131$		& Odd, standard 8-form	& $=8(16)+3$		\\
				&						&							& $=8(3*5+1)+3$	\\
\end{tabular}
\end{table}

In this cycle the base of the number decreases from 16 to 8, the index increases from 5 to 16, and the value of the number increases from 87 to 131.
 
\newpage 
\subsection*{Appendix V - Analysis of General 16-Cascade}
\begin{table}[!htbp]
\centering 										
\begin{tabular}{lllp{5cm}} 						
Step	& Form									& Resulting form			& Resulting parity \\
\hline
0		& $16n + 7$ 							& $= 16n + 7$				& Odd, standard 16-form \\
1		& $3(16n + 7) + 1= 48n + 22$		&							& Even \\
2		& $(48n + 22) / 2= 24n + 11$ 		& $= 8(3n + 1) + 3$		& Odd, standard 8-form \\
3		& $3(24n + 11) + 1= 72n + 34$									& Even \\
4		& $(72n + 34) / 2= 36n + 17$ 		& $= 4(9n + 4) + 1$		& Odd, standard 4-form Note: $9n + 4 = 3(3n + 1) + 1$ \\
5		& $3(36n + 17) + 1= 108n + 52$	&							& Even \\
6		& $(108n + 52) / 2 = 54n + 26$ 		& $=2(27n + 13)$			& Even, standard 2-form Note: $27n + 13 = 3(9n + 4) + 1$ \\
7		& $(54n + 26) / 2 = 27n + 13$ 		& $= 27n + 13$			& Mixed form, Odd or even depending on n. \\
\end{tabular}
\end{table}

\end{document}